\RequirePackage{fix-cm}
\documentclass[smallextended]{svjour3}       
\smartqed  
\usepackage{graphicx}
\usepackage{bibentry}
\usepackage{graphics}
\usepackage{url}
\usepackage{amsmath,amssymb}
\usepackage{color}

\newcommand{\RR}{\mathbb{R} }
\newcommand{\LL}{\mathbb{L} }

\newcommand{\St}{\mathbb{S} }
\newcommand{\Wt}{\mathbb{W} }

\newcommand{\isun}{i su(N)}

\newcommand{\beq}{\begin{equation}}
\newcommand{\eeq}{\end{equation}}
\newcommand{\beqn}{\begin{eqnarray}}
\newcommand{\eeqn}{\end{eqnarray}}
\newtheorem{thm}{Theorem}
\newtheorem{rem}{Remark}

\newcommand{\e}{\epsilon}
\newcommand{\hmu}{\hat{\mu}}
\newcommand{\bproof}{{\noindent {\bf Proof} }}
\newcommand{\eproof}{{\ensuremath{\square}}}

\journalname{Acta Applicandae Mathematicae}
\begin{document}

\title{Critical points of the optimal quantum control landscape: a propagator 
approach
}


\author{Tak-San Ho \and Herschel Rabitz \and Gabriel Turinici}


\institute{
T.-S. Ho \at
Department of Chemistry, Princeton University,
Princeton, NJ 08544, USA \\
\email{tsho@princeton.edu}  
\and
H. Rabitz \at
Department of Chemistry, Princeton University,
Princeton, NJ 08544, USA \\
\email{hrabitz@princeton.edu}    
\and
G. Turinici \at
CEREMADE, Universit\'e Paris Dauphine, 
Place du Marechal de Lattre de Tassigny, 75016 Paris, France \\
\email{Gabriel.Turinici@dauphine.fr}           
}

\date{Received: date / Accepted: date}

\maketitle

\begin{abstract}
Numerical and experimental realizations of  quantum control 
are closely connected to the 
 properties of the mapping from the control to the unitary propagator~\cite{science04rabitz}. 
For bilinear quantum control problems, no general results are available to fully determine
when this mapping is singular or not. In this paper we give sufficient conditions, 
in terms of elements of the evolution semigroup, for a 
trajectory to be non-singular. We 
identify two lists of ``way-points" that, when reached, ensure the non-singularity of
the control trajectory. It is found that under appropriate hypotheses one of those lists does not depend 
on the values of the coupling operator matrix. 
\keywords{quantum control \and singular control \and landscape analysis in quantum control}
 \PACS{32.80.Qk \and 02.30.Yy}
 \subclass{15A57 \and 49J20 \and 65Z05}
\end{abstract}

\section{Introduction}

Manipulating the evolution of physical systems at the quantum level
has been a longstanding goal from the very beginnings of laser technology.
Investigations in this area accelerated after the introductions of
(optimal) control theory tools~\cite{hbref110}, which greatly contributed to the
first positive experimental
results, see~\cite{hbref1,hbref3,hbref4,hbref5,hbref6,hbref112,hbref113} and references herein. 

Consider a quantum system in the presence of a control field $\e(t) \in \RR, \ t \ge 0$. In the electric dipole approximation 
and  in the  density matrix formulation, the underlying time-dependent density
matrix $\rho(t)$ satisfies {the Liouville-von Neumann equation}
\beqn \label{eqn:densityevolution}
& \ &
i \frac{\partial}{\partial t} \rho(t) =
[  H_0  -\e(t)\mu, \rho(t)]
\\ & \ &
\rho(t=0) = \rho_0
\nonumber
\eeqn
where $H_0$ is the field-free Hamiltonian
(including the potential) and $\mu$ is the dipole moment of the system.
In the analysis, we will suppose that the system  contains  $N$ levels,  thus,
 $H_0$, $\mu$ and $\rho(t)$ are all $N\times N$ Hermitian matrices.  In the following, for simplicity,
we further assume that $H_0$ and $\mu$ are real and symmetric matrices.
Moreover, to avoid trivial settings, the dipole moment operator $\mu$ is assumed to have zero trace, i.e.,
\beq \label{eq:muzero}
Tr(\mu) =0.
\eeq
{In the following we adopt the notation $\isun$ to denote the set of all complex $N \times N$ Hermitian matrices with zero trace.}
Recall that $\rho(t)= U(t,0) \rho(0) U^*(t,0)$, where
the unitary propagator $U(t,0) \in U(N)$, $t\ge0$ is the solution of {the time-dependent Sch{\"o}dinger equation}
\beqn \label{eqn:propagator}
& \ &
i \frac{\partial}{\partial t} U(t,0) =
\Big(H_0  -\e(t)\mu \Big)U(t,0)
\\ & \ &
U(t=0,0) = I_N.
\nonumber
\eeqn
The control goal can be expressed in terms of a (self-adjoint) operator $O$ in that the 
corresponding functional to be optimized (maximized) is  
$\langle O\rangle(T) = {\rm Tr}\Big(\rho(T)O\Big)$.  Accordingly, a maximum value of 
$\langle O\rangle(T)$ indicates that the control is of optimal quality and {\it vice-versa}.

At the heart of a good understanding of any quantum control
problem, including the efficiency, accuracy, 
and stability of quantum control experiments as well as simulations,
lies essential features of the control landscape, depicting
the functional dependence of the observable on the control field,
i.e.  the mapping: $\e(t) \mapsto \langle O\rangle(T)$. The 
slopes (gradient) and the curvatures (Hessian) are key characteristics of the
control landscape. In particular, the topologies at and around the critical
points, in which the slopes vanish, can shed important light on the questions,
for example, of what makes quantum control experiments (simulations) 
apparently ``easy to perform'' and why quantum control beats the ``curse of 
dimensionality'' (i.e., overcomes the anticipated exponentially growing effort 
required when searching over increasing numbers of control variables)~\cite{JPPAhr2006}. A
thorough study of the underlying control landscape should give insight into 
the classes of future feasible quantum control experiments, especially 
those involving complex molecules and materials.

\section{Motivation: landscape analysis and beyond}

The search procedures optimize the control quality 
 $\langle O\rangle(T)$ with respect to variations in the control
 $\e(t)$. The mapping 
 $\e(t) \mapsto \langle O\rangle(T)$ is the composition of two maps: a control
map 
 $\e(t) \mapsto U(T,0)$ and a "kinematic" map
 $U(T,0) \mapsto \langle O\rangle(T)$.
A key issue is whether there exist critical points that are suboptimal solutions. To analyse the critical points it is expedient to consider
 these two maps separately.

Previous analyses~\cite{science04rabitz} have investigated
the mapping from the propagator to the observable i.e. $ U(T,0) \mapsto \langle O\rangle(T)$; 
the critical points correspond to the equation 
\begin{equation}
\frac{\delta \langle O\rangle(T)}{\delta U(T,0)}=0.
\end{equation}
Assuming controllability hypothesis~\cite{altafini,controllability1}, this
equation is satisfied if and only if
 \begin{equation}
 [U^*(T,0)OU(T,0), \rho(0)]=0,
 \end{equation}
 i.e, we obtain a necessary and sufficient condition for the kinematic critical points~\cite{PRAgslk1998}. The optimization "lanscape", when looked upon in terms of this so called "kinematic" mapping, is therefore exempt
of suboptimal critical points.

To extend the landscape analysis  
beyond the kinematic setting, i.e. to investigate the dependence on the control itself,
we need to analyse the solutions of the following critical point condition~\cite{JPPAhr2006}
\begin{equation}
\frac{\delta \langle O\rangle(T)}{\delta\epsilon(t)}=0.
\end{equation}
The solutions of this equation are among the critical points 
of the control map $\e(t) \mapsto U(T,0)$; it was shown in~\cite{2009arXiv0907.2354}
that this mapping is nonsingular {when all matrix entries of the dipole moment
\beq \label{eq:muhat}
\hmu(t) =  U^*(t,0) \mu U(t,0),
\eeq
 as functions of
$t\in[0,T]$, are linearly independent
(to the extent possible), or equivalently
if the dipole moment matrices $\hmu(t)$'s
span the space of all zero-traced Hermitian matrices $\isun$.}

In previous works, it has
always been assumed that ``full controllability'' implies ``fully linear 
independence'' of the dipole moment matrix entries. An assertion of
the linear independence of dipole moment matrix entries is also essential 
for understanding the landscape of unitary transformation control 
problem~\cite{hdr2007ut}.
In this paper, we present a rigorous analysis of this assertion and identify properties of 
the trajectory, in terms of the propagators, that are sufficient to ensure the
linear independence of entries of $\hmu(t)$. 
{In addition,
the analysis of which controls correspond to a non-singular mapping $\e(t) \mapsto U(T,0)$
or $\e(t) \mapsto \langle O\rangle(T)$ is also relevant, e.g., for stabilization purposes, error 
cancellation, maintaining coherence, etc.}
\section{Dipole dependent way-points}

Denote by $\LL$ the Lie algebra generated by $-iH_0$ and $-i\mu$ as a sub-algebra of
$u(N)$.  We will suppose that  {the system is density matrix controllable, or equivalently
\beq \label{hyp:cont}
\LL = su(N)\ \mbox{or}\ \LL = u(N).
\eeq
}Note that when 
the system is density matrix controllable, then
it is also wave function controllable, cf~\cite{alessandro1,alessandrobook}.

\begin{thm} 
Under assumptions~\eqref{eq:muzero} and \eqref{hyp:cont}
there exists a set $\{U_1,...,U_{2N^2-2N}\}$ of unitary propagators (that we will call  
the set of "way-points")
such that if 
$U(t,0)$ visits all of them i.e.
$U(t_k,0) = U_k$ for some times $t_k$, $k=1,...,{2N^2-2N}$ the components
$f_{ij}(t)$ of the matrix $\hmu(t)$ in~\eqref{eq:muhat} are linearly independent as
functions of time.

Consequently, under the hypotheses~\eqref{eq:muzero} and \eqref{hyp:cont} for $T$ sufficiently large a control field 
$\e(t)$, $t\in [0,T]$ exists
such that functions $f_{ij}(t)$ are linearly independent (or equivalently such that
$\hmu(t)$ spans the set of all zero-traced Hermitian matrices).
\end{thm}

\begin{rem}
Note that here there is no information required on the structure of the matrices $H_0$ or $\mu$.
\end{rem}
\bproof
Since the Hermitian matrix $\mu$ has zero trace 
it cannot be a constant and thus it has at least two different eigenvalues
$\lambda_1$ and $\lambda_2$.

Let us begin by some general conventions: for any $i \neq j$  and any matrix $M$ 
we denote by
$M_{<i,j>}$ its $2\times2$ sub-matrix formed by the $i$-th and $j$-th rows and columns:
\beq
M_{<i,j>} = \begin{pmatrix} 
M_{ii} & M_{ij} \\
M_{ji} & M_{jj}
\end{pmatrix}.
\eeq
Moreover,
when we alter the $N\times N$ identity matrix by putting 
in the $i$-th and $j$-th rows and columns a given $2\times2$ matrix $D$,
we will denote the resultant $N\times N$ matrix by $D^{<i,j>}$~:

\beq
D^{<i,j>}=  
\begin{matrix}
\begin{matrix}  \ & \ \  & i & \ \ & j \ \ & \ &  \end{matrix}
\\
\left( 
\begin{matrix} 1 \\ 0 \\ \ \\0 \\ \ \\ 0 \\ \ \\ 0 \end{matrix}
\ \ 
\begin{matrix} 0 \\ 1 \\ \ \\0 \\ \ \\ 0 \\ \ \\ 0 \end{matrix}
\ \ 
\ \ 
\begin{matrix} 0 \\ 0 \\ \ \\D_{11} \\ \ \\ D_{21} \\ \ \\ 0 \end{matrix}
\ \ 
\ \ 
\begin{matrix} 0 \\ 0 \\ \ \\D_{12} \\ \ \\ D_{22} \\ \ \\ 0 \end{matrix}
\ \ 
\ \ 
\begin{matrix} 0 \\ 0 \\ \ \\0 \\ \ \\ 0 \\ \ \\ 1 \end{matrix}
\right)
\ \ 
\begin{matrix} \ \\ \ \\ \ \\i \\ \ \\ j \\ \ \\ \ \end{matrix}
\end{matrix}
\eeq

We choose $U_k$ in the following way: let $l=1,...,\frac{N(N-1)}{2}$ index the couples $(i,j)$,
$i,j=1,...,N$ with $i < j$.

- $U_{4l+1}$ is such that $\Big(\hat{\mu}(t_{4l+1})\Big)_{<i,j>} = 
 \left(
\begin{array}{cc} 	
\lambda_1 & 0 \\
0         &\lambda_2
\end{array}  
\right) 
$;

- $U_{4l+2}= U_{4l+1}  \left(
\begin{array}{cc} 	
0 & 1 \\
1 & 0
\end{array} 
\right)^{<i,j>}
$. Then $\Big(\hat{\mu}(t_{4l+2})\Big)_{<i,j>} = 
 \left(
\begin{array}{cc} 	
\lambda_2 & 0 \\ 
0         &\lambda_1
\end{array} 
\right) 
$ and all other entries are identical to those of 
$\hat{\mu}(t_{4l+1})$.

- $U_{4l+3}=  \frac{1}{2} U_{4l+1} \left(
\begin{array}{cc} 	
1 & 1 \\
-1 & 1
\end{array} 
\right)^{<i,j>}
$. Then 
\beq
\Big(\hat{\mu}(t_{4l+3})\Big)_{<i,j>} = 
\frac{1}{2}  \left(
\begin{array}{cc} 	
\lambda_1+\lambda_2 & \lambda_1-\lambda_2 \\
\lambda_1-\lambda_2  &\lambda_1+\lambda_2
\end{array} 
\right) 
\eeq
and all other entries are identical to those of 
$\hat{\mu}(t_{4l+1})$.

- $U_{4l+4}=  \frac{1}{2} U_{4l+1} \left(
\begin{array}{cc} 	
i & 1 \\
1 & i
\end{array} 
\right)^{<i,j>}
$. Then 
\beq
\Big(\hat{\mu}(t_{4l+4})\Big)_{<i,j>} = 
\frac{1}{2}  
\left(
\begin{array}{cc} 	
\lambda_1+\lambda_2 & (\lambda_2-\lambda_1)i \\
(\lambda_1-\lambda_2)i  &\lambda_1+\lambda_2
\end{array} 
\right) 
\eeq
and all other entries are identical to those of 
$\hat{\mu}(t_{4l+1})$.

Note that all $\hmu(t_k)= U_k^* \mu U_k$ are Hermitian (and of zero trace).
Let us take a 
zero trace Hermitian matrix $Z$ 
such that 
\beq
\langle Z, \hat{\mu}(t_k) \rangle = 0, \ \forall k=1,...,4N.
\eeq
where we used the canonical
scalar product of Hermitian matrices $\langle A, B\rangle = Tr(A^*B)=Tr(AB)$.
We recall also the definition of the norm induced by 
this scalar product $\| A \|= \sqrt{\langle A, A\rangle}$.

We denote $
\left(
\begin{array}{cc} 	
x & y \\
y^*  &z
\end{array} 
\right) =Z_{<i,j>}
$; since $\langle Z, \hat{\mu}(t_k) \rangle = 0$ for $k=4l+1,...,4l+4$ it follows that
$\langle Z_{<i,j>}, (\hat{\mu}(t_k))_{<i,j>} \rangle$ is the same for $k=4l+1,...,4l+4$. The consequence is that 
$\lambda_1 x + \lambda_2 z = \lambda_2  x + \lambda_1 z $ which together
with $\lambda_1\neq \lambda_2$ implies $x=z$. Moreover we also have
$\lambda_1 x + \lambda_2 x = \frac{1}{2}\Big(
2(\lambda_1  + \lambda_2 )x + (\lambda_1  - \lambda_2 )(y +y^*)
\Big)
$ thus $y +y^*=0$. Additionally, we have
$\lambda_1 x + \lambda_2 x = \frac{1}{2}\Big(
2(\lambda_1  + \lambda_2 )x + (\lambda_1  - \lambda_2 )i(y^*-y)
\Big)
$, thus $y =y^*$.  From $y +y^*=0$ and $y =y^*$ we infer that $y=0$ and
 obtain $Z_{<i,j>} = x  I_2$ ($I_2$ is the $2\times2$ unit matrix).
Since this is true for arbitrary $i<j$ we obtain $Z = x I_N$. But, since $Tr(Z)=0$ it
follows $Z=0$, q.e.d.

For the second part of the conclusion we use the fact that under the hypotheses~\eqref{eq:muzero} and \eqref{hyp:cont} 
the system is controllable thus 
for $T$ sufficiently large a control field 
$\e(t)$, $t\in [0,T]$ exists
such that functions $f_{ij}(t)$ are linearly independent.

\section{Dipole independent way-points}

The result of the previous section can be interpreted as follows: as long as the control field 
ensures that the propagator will visit some specific unitary transformations, then the required linear independence
of the time-dependent elements of $\hmu(t)$ is satisfied. 
Note that this set of unitary transformations (i.e., propagators) explicitly depends on the dipole moment.
In the laboratory, procedures have been designed that are able to control the quantum evolution even in the absence of precise information
on the dipole~\cite{hbref110}. It is therefore possible in principle 
to experimentally find a control passing through a specified list of propagators when this list
does not depend on the dipole entries. A universal procedure can therefore be implemented that will
ensure the non-singularity of the mapping $\e(t) \mapsto U(T,0)$ even when the dipole is unknown.
Beyond the question of the landscape analysis, such a procedure can be useful in additional
circumstances when the mapping is required to be non-singular, e.g. for stabilization purposes, etc.
We will therefore investigate the following circumstance: suppose that the system is
controllable and thus one can experimentally implement  arbitrary propagators $U\in SU(N)$. Can we find a list
of propagators (as in the previous section) which are independent of the precise values of entries of $\mu$ ?

\begin{thm} 
Under assumptions~\eqref{eq:muzero} and \eqref{hyp:cont} 
there exists a set of "way-points"
$\Wt \subset SU(N)$ independent of $\mu$
such that if 
$U(t,0)$ visits all propagators in $\Wt$ (i.e.
for all $U_k \in \Wt$ there exists $t_k$ with
$U(t_k,0) = U_k$)  the components
$f_{ij}(t)$ of the matrix $\hmu(t)$ in~\eqref{eq:muhat} are linearly independent as
functions of time.
\end{thm}
\bproof
We denote for any $U \in SU(N)$:
\beq
{\mathcal C}_{U}=\{ (Z,\mu)\in (\isun)^2 \ | \ \|Z\|=\|\mu\|=1, \ Tr(Z U^* \mu U)\neq 0 \}.
\eeq
This set is open because its complement is the solution of a linear equation in the entries of the matrices $Z$ and $\mu$. 

We prove now that for any $Z,\mu \in \isun$ of unit norm there exists
a $U\in SU(N)$ such that $Tr(Z U^* \mu U)\neq 0$. To see this, diagonalize $Z$ and $\mu$: 
$Z = U_1^* D_1 U_1$ and $\mu = U_2^* D_2 U_2$ with $U_1,U_2$ unitary and $D_1, D_2$ diagonal; we will denote by $d^a_k$ the k-th diagonal entry of $D^a$. 
Recall that since $Z,\mu$ are zero-traced $\sum_{k=1 }^N d^a_k = 0$ for $a=1,2$. Therefore, in each set $(d^a_k)_{k=1}^N$ there are some elements that are strictly negative and some strictly positive. 
We will also suppose that the diagonalization is done in such a way that the diagonal elements are ordered from the lowest (that is strictly negative) to the highest, which is strictly positive.
Then,
\beq
Tr(Z U^* \mu U) = Tr( U_1^* D_1 U_1 U^* U_2^* D_2 U_2 U) = Tr(D_1 V^* D_2 V)
\eeq
with $V =U_2 U U_1^* \in SU(N)$.   Take now $U$ such that $V$ is a permutation matrix corresponding to permutation $\sigma$. Then  
\beq
 Tr(D_1 V^* D_2 V)= \sum_{k=1}^N d^1_k d^2_{\sigma(k)}.
\eeq
If $\sum_{k=1}^N d^1_k d^2_{\sigma(k)}=0$ for any permutation $\sigma$ then
it follows in particular that
$\sum_{k=1}^N d^1_k d^2_{k} = \sum_{k=1}^N d^1_k d^2_{N-k}$, but the first member is
always superior to the second with the equality implying that one of the vectors  $d^a_k$ is constant  which leads to a contradiction with the zero-trace hypothesis.

We therefore obtain that $\cup_{U\in SU(N)} {\mathcal C}_{U}$ is a covering of the compact set 
$\{ (Z,\mu)\in (\isun)^2 \ | \ \|Z\|=\|\mu\|=1\}$ and thus one can extract a finite covering 
i.e. a finite set $\Wt$ such that
\beq
\cup_{U \in \Wt} {\mathcal C}_{U} =\{ (Z,\mu)\in (\isun)^2 \ | \ \|Z\|=\|\mu\|=1\}.
\eeq 
This means that for any $Z$ and $\mu$ in $\isun$ (note that there is no need to impose the norm condition) there exists $U \in \Wt$ such that 
$Tr(Z U^* \mu U)\neq 0$. Expressed otherwise, for any $\mu \in \isun$ 
there does not exist $Z \in \isun$ such that 
$Tr(Z U^* \mu U)=0$ for any $U \in \Wt$. This implies that
$\{ U^* \mu U \ | \ U \in \Wt \}$ form an independent set of vectors in $\isun$, q.e.d.
\eproof

\begin{rem}
Theorem 2 only assures that such a set $\Wt$ exists but does not give yet precise information on its size (which can be large) nor does it provide a constructive approach to identify in practice 
$\Wt$. The following result addresses these questions upon imposing an additional constraint on the dipole $\mu$.
\end{rem}

\begin{thm} 
Under assumptions~\eqref{eq:muzero} and \eqref{hyp:cont}
and
\beq \label{hyp:munonzero}
\mu_{ij}\neq 0, \ \forall i \neq j.
\eeq
there exists a set of "way-points"
$\St \subset SU(N)$ 
such that if 
$U(t,0)$ visits all propagators in $\St$ (i.e.
for all $U_k \in \St$ there exists $t_k$ with
$U(t_k,0) = U_k$)  the components
$f_{ij}(t)$ of the matrix $\hmu(t)$ in~\eqref{eq:muhat} are linearly independent as
functions of time. Moreover the set $\St$ can be chosen
to be the same for all coupling operators $\mu$ satisfying the 
hypotheses~\eqref{eq:muzero}, \eqref{hyp:cont}, \eqref{hyp:munonzero}.
\end{thm}
\bproof 
Suppose that the entries of $\hmu(t)$ are not linearly independent 
i.e. a matrix $Z$, $Tr(Z)=0$ exists such that $Z$ is orthogonal to all
matrices $U^* \mu U$, $\forall U \in \St$, i.e., 
$Tr(Z U^* \mu U) = 0$ $\forall U \in \St$.

Let us compute  $U^* \mu U$:
for $\mu=\left( 
\begin{matrix}  
\mu_{\bf 11} & \mu_{\bf 12} \\
\mu_{\bf 21}   & \mu_{\bf 22}\\
\end{matrix}
\right)
$ 
(here $\mu_{\bf 11}$ is a $2 \times 2$ matrix, 
$\mu_{\bf 22}$ a $N-2 \times N-2$ matrix, etc.)
and $U= W^{<1,2>}=
\left( 
\begin{matrix}  
W & 0 \\
0   & I_{N-2} &\\
\end{matrix}
\right)
$
\beq
U^* \mu U = 
\left( 
\begin{matrix}  
W^*\mu_{\bf 11}W & W^*\mu_{\bf 12} \\
\mu_{\bf 21}W   & \mu_{\bf 22}\\
\end{matrix}
\right).
\eeq
A typical example of $W$ is 
$\left( 
\begin{matrix}  
\cos(\phi) & \sin(\phi) \\
\sin(\phi)   & -\cos(\phi)\\
\end{matrix}
\right).
$
The entries of $W^*\mu_{\bf 11}W$ will contain terms of
second order in $\cos(\phi)$ and $\sin(\phi)$:
$\cos^2(\phi)= \frac{\cos(2\phi)+1}{2}$ and 
$\sin^2(\phi)=\frac{1-\cos(2\phi)}{2}$, 
$\sin(\phi)\cos(\phi)=\frac{\sin(2\phi)}{2}$; the entries of 
$W^*\mu_{\bf 12}$ and $\mu_{\bf 21}W$ will contain only
first order terms in 
$\cos(\phi)$ and $\sin(\phi)$ while $\mu_{\bf 22}$ does not
depend on $\phi$ at all. 

We recall that for any $\theta \in \RR$ the functions
$1$, $\cos(\theta)$,
 $\sin(\theta)$, 
 $\cos(2\theta)$, 
 $\sin(2\theta)$ are linearly independent, which will be used in the following 
fashion
\begin{lemma} \label{lemma:theta}
There exist five real values of $\theta$,
for instance $\Theta=\{0, \pi/3,\pi/2, \pi, 3\pi/3 \}$
such that 
\beqn
& \ & \alpha_1 + \alpha_2 \cos(\theta) +  \alpha_3 \sin(\theta) \nonumber\\
& \ & 
+  \alpha_4 \cos(2\theta)
+ \alpha_5 \sin(2\theta)=0
\ \forall \theta \in \Theta
\eeqn  
implies 
\beq
\alpha_1=...=\alpha_5=0.
\eeq
\end{lemma}

We now introduce into the set $\St$ the matrices 
\beq
U_{\theta,i;j}=\left(
\begin{array}{cc} 	
0 & \cos(\theta)+i\sin(\theta) \\
\cos(\theta)-i\sin(\theta) &0
\end{array} 
\right)^{<i,j>},
\eeq
 for any $i\neq j$ and $\theta \in \Theta$.

Since $Tr(Z U_{\theta,i;j}^* \mu U_{\theta,i;j}) = 0$ $\forall \theta \in \Theta$
by Lemma~\ref{lemma:theta}, the coefficients of $\sin{2\theta}$
and $\cos{2\theta}$ in $Tr(Z U^* \mu U) = 0$ will be zero.
The coefficient of $\sin{2\theta}$ turns out to be $\mu_{ij}(Z_{ij}+Z_{ji})$
and that of $\cos{2\theta}$ is $i \mu_{ij}(Z_{ij}-Z_{ji})$. It follows from~\eqref{hyp:munonzero}
that $Z_{ij}=Z_{ji}=0$.

Let us add to the set $\St$ the matrices 
\beq
V_{\theta,i;j}=\left(
\begin{array}{cc} 	
\cos(\theta) & \sin(\theta) \\
\sin(\theta) & -\cos(\theta) 
\end{array} 
\right)^{<i,j>}, 
\eeq
 for any $i=1,...,N-1$, $j=i+1$ and $\theta \in \Theta$.
We use again
$Tr(Z V_{\theta,i;j}^* \mu V_{\theta,i;j}) = 0$ $\forall \theta \in \Theta$ and identify the coefficients
of $\sin{2\theta}$
and $\cos{2\theta}$ which are respectively $\mu_{ij}(Z_{ii}-Z_{jj}) +\frac{\mu_{ii}-\mu_{jj}}{2}(Z_{ij}+Z_{ji}) $
and  $-\mu_{ij}(Z_{ij}+Z_{ji}) +\frac{\mu_{ii}-\mu_{jj}}{2} (Z_{ii}-Z_{jj})$. We obtain $Z_{ii}=Z_{jj}$.
Since this is true for any $i=1,...,N-1$ from $Tr(Z)=0$ it follows that $Z_{ii}=0$ $\forall i\le N$.

\begin{rem}
Although the hypothesis~\eqref{hyp:munonzero} may seem strong we believe
that it is rather a technical requirement which could hopefully be relaxed 
in the future.
\end{rem}

\section{Conclusion and discussion}

Following previous works in~\cite{science04rabitz} 
which study the properties of the control input-output map we investigate in this paper
the singularity of the control to propagator map through the study of
the linear independence of entries of the time-dependent coupling operator
$\hmu(t)$ (cf.~\eqref{eq:muhat}).  We provide several
sufficient conditions in terms of the evolution semigroup points that are reached during
the propagation of Eq.~\eqref{eqn:densityevolution}. We expect that the criterion in \eqref{hyp:munonzero} can be relaxed to accommodate even more general coupling operators.


\section*{Acknowledgements}

GT acknowledges support from the ANR-06-BLAN-0052 program 
and the MicMac project of INRIA Rocquencourt. TSH and HR were supported by U.S. Department of Energy.

\bibliographystyle{spmpsci}      
\bibliography{refs}   

\end{document}